\documentclass[11pt,american,english]{article}

\usepackage[latin9]{inputenc}
\usepackage[a4paper]{geometry}
\usepackage{amsmath}
\usepackage{amssymb}
\usepackage{setspace}
\makeatletter

\providecommand{\tabularnewline}{\\}

\@ifundefined{date}{}{\date{}}
\usepackage{array}
\usepackage{mathrsfs}
\usepackage{color}

\makeatother

\usepackage{babel}
\begin{document}

\title{\textbf{Conjunctions of Three ``Euler Constants'' in}\\
\textbf{Poisson-Related Expressions }}

\author{Michael R. Powers\thanks{Department of Finance, School of Economics and Management, and Schwarzman
College, Tsinghua University, Beijing, China 100084; email: powers@sem.tsinghua.edu.cn.}}

\date{July 1, 2025}
\maketitle
\begin{abstract}
\begin{singlespace}
\noindent Three mathematical constants bear the name of the venerable
Leonhard Euler: Euler's number, $e=2.718281\ldots$; the Euler-Mascheroni
constant, $\gamma=0.577216\ldots$; and the Euler-Gompertz constant,
$\delta=0.596347\ldots$. In the present work, we consider two joint
appearances of these constants, one in a well-known equation of Hardy
(interpretable in connection with inverse second moments of the Poisson
probability distribution), and the other from a sequence of probabilities
generated by recursively conditional Exponential (i.e., Poisson-event
waiting-time) distributions. In both cases, we explore generalizations
of the initial observations to offer more comprehensive results, including
extensions of Hardy's equation.\medskip{}

\end{singlespace}
\end{abstract}

\section{Introduction}

\noindent The influence of Leonhard Euler (1707-1783) is pervasive
in modern mathematics. Through a prodigious body of research (886
published papers), he not only contributed profoundly to the research
of his time, but also laid foundational ideas for new subfields (e.g.,
topology) that developed after his passing. Nowhere was Euler's impact
greater than in the early study of special functions and analytic
number theory. Among the many results and objects bearing his name
today are three mathematical constants. In decreasing order of familiarity,
these are: Euler's number, $e=2.718281\ldots$; the Euler-Mascheroni
constant, $\gamma=0.577216\ldots$; and the Euler-Gompertz constant,
$\delta=0.596347\ldots$.

Euler's number, like the constant $\pi$, is ubiquitous. Although
$\pi$ already was known to the ancients, $e$ was not defined formally
until Jacob Bernoulli used it to analyze compound interest in the
late 17th century. The tradition of naming this quantity after Euler
appears to have come from his extensive study of the constant as the
base of the exponential function as well as general acceptance of
his use of the symbol $e$.

The Euler-Mascheroni constant was discovered by Euler, but shares
naming honors with Lorenzo Mascheroni because of the latter's detailed
approximation of its value. Euler first determined the constant in
his investigation of the harmonic series, where it appears as
\begin{equation}
\gamma=\underset{n\rightarrow\infty}{\lim}\left({\displaystyle \sum_{i=1}^{n}}\dfrac{1}{i}-\ln\left(n\right)\right).
\end{equation}
In addition to connections with various special functions (gamma,
exponential-integral, zeta, Bessel, etc.), it may be expressed in
numerous integral forms, including 
\begin{equation}
\gamma=-{\displaystyle \int_{0}^{\infty}e^{-t}\ln\left(t\right)dt}
\end{equation}
and
\[
\gamma={\displaystyle \int_{0}^{1}\left(\dfrac{1}{\ln\left(t\right)}+\dfrac{1}{1-t}\right)dt}.
\]

Finally, the Euler-Gompertz constant also was discovered by Euler,
who identified it as the (regularized) sum of the divergent Wallis
series,
\[
\delta=0!-1!+2!-3!+\cdots
\]
(famously included in the first letter of Srinivasa Ramanujan to G.
H. Hardy; see Berndt and Rankin, 1995). This constant is named jointly
with Benjamin Gompertz because of its appearance in calculations related
to Gompertz's work in mortality/survival analysis, although it appears
Gompertz never investigated the quantity closely (see Lagarias, 2013).
Often, it is rendered as
\[
\delta={\displaystyle \int_{0}^{\infty}e^{-t}\ln\left(t+1\right)dt},
\]
an integral similar to (2). Other expressions include
\[
\delta={\displaystyle \int_{0}^{\infty}\dfrac{e^{-t}}{t+1}dt},
\]
\[
\delta={\displaystyle \int_{0}^{1}\dfrac{1}{1-\ln\left(t\right)}dt},
\]
and
\[
\delta=-e\textrm{Ei}\left(-1\right),
\]
where $\textrm{Ei}\left(u\right)={\textstyle \int_{-\infty}^{u}t^{-1}e^{t}}dt$
denotes the exponential-integral function.

Although the irrationality of $e$ was demonstrated by Euler, both
$\delta$ and $\gamma$ so far have resisted attempts to prove that
they too are irrational. Consequently, there is considerable research
interest in investigating their analytic properties. In recent years,
Aptekarev (2009) proved that at least one of the two constants must
be irrational, and Rivoal (2012) showed further that at least one
of them must be transcendental.

In the present work, we consider two intriguing conjunctions (i.e.,
joint appearances) of the three ``Euler constants'' in expressions
associated with the Poisson probability distribution. The first case,
addressed in Section 2, involves a well-known equation of Hardy (1949)
that may be interpreted in relation to inverse second moments of the
Poisson random variable. The second example, presented in Section
3, arises from a sequence of probabilities generated by recursively
conditional Exponential (i.e., Poisson-event waiting-time) distributions.
In both cases, we explore generalizations of the initial observations
to place them within a broader context from which more comprehensive
results, including extensions of Hardy's equation, may be derived.

\section{Hardy's Equation}

\noindent Consider the power-series expansion of $\textrm{Ei}\left(u\right)$,
given by
\begin{equation}
\textrm{Ei}\left(u\right)=\gamma+\ln\left|u\right|+{\displaystyle \sum_{i=1}^{\infty}\dfrac{u^{i}}{i\cdot i!}},
\end{equation}
for $y\in\mathbb{R}\setminus\left\{ 0\right\} $. Evaluating this
expression at $u=-1$, it is easy to see that
\[
\textrm{Ei}\left(-1\right)=\gamma+{\displaystyle \sum_{i=1}^{\infty}\dfrac{\left(-1\right)^{i}}{i\cdot i!}}
\]
\[
\Longleftrightarrow\delta=-e\left(\gamma-{\displaystyle \sum_{i=1}^{\infty}\dfrac{\left(-1\right)^{i+1}}{i\cdot i!}}\right),
\]
a result derived by Hardy (1949) and well known for the conjunction
of all three Euler constants. This equation may be rearranged more
elegantly as
\begin{equation}
{\displaystyle \sum_{i=1}^{\infty}\dfrac{\left(-1\right)^{i+1}}{i\cdot i!}}=\dfrac{\delta}{e}+\gamma,
\end{equation}
which enjoys the alternative integral form,
\begin{equation}
{\displaystyle \int_{0}^{1}}\dfrac{1-e^{-t}}{t}dt=\dfrac{\delta}{e}+\gamma.
\end{equation}

One shortcoming of expressions (4) and (5) is that the left-hand sides
do not possess simple, intuitive interpretations. To improve matters,
one can observe that
\[
{\displaystyle \sum_{i=1}^{\infty}\dfrac{\left(-1\right)^{i+1}}{i\cdot i!}}=e{\displaystyle \sum_{i=0}^{\infty}}\dfrac{e^{-1}1^{i}}{i!}\dfrac{\left(-1\right)^{i}}{\left(i+1\right)^{2}}
\]
\[
=e\textrm{E}_{N\mid\lambda=1}\left[\dfrac{\left(-1\right)^{N}}{\left(N+1\right)^{2}}\right],
\]
where $N\mid\lambda\sim\textrm{Poisson}\left(\lambda=1\right)$. Setting
\begin{equation}
\mu_{u,k}=\textrm{E}_{N\mid\lambda=1}\left[\dfrac{u^{N}}{\left(N+k\right)^{2}}\right]
\end{equation}
then allows us to rewrite (4) as
\[
e\mu_{-1,1}=\dfrac{\delta}{e}+\gamma
\]
\begin{equation}
\Longleftrightarrow\mu_{-1,1}=\dfrac{\delta+e\gamma}{e^{2}}.
\end{equation}

Introducing the quantity in (6) also facilitates generalizations of
Hardy's equation through the parameters $u$ and $k$. For $u=-1$,
(6) is essentially the alternating inverse second moment of $N+k$,
where $N\mid\lambda\sim\textrm{Poisson}\left(\lambda=1\right)$, and
is finite only if $k>0$. We consider all $k\in\mathbb{\mathbb{Z}}^{+}$
for $u=-1$ in Subsection 2.1, and then address the case of $u=1$
in Subsection 2.2.

\subsection{General $\boldsymbol{k}$ for $\boldsymbol{u=-1}$}

\noindent To evaluate $\mu_{-1,k}$ for $k\in\left\{ 2,3,\ldots\right\} $,
we begin with $k=2$ and proceed iteratively. From the definition
in (6), we have

\noindent 
\[
\mu_{-1,2}={\displaystyle \sum_{i=0}^{\infty}}\dfrac{e^{-1}1^{i}}{i!}\dfrac{\left(-1\right)^{i}}{\left(i+2\right)^{2}}
\]
\[
={\displaystyle \sum_{\ell=1}^{\infty}}\dfrac{e^{-1}1^{\ell-1}}{\ell!}\dfrac{\left(-1\right)^{\ell-1}\ell}{\left(\ell+1\right)^{2}}
\]
\[
=-{\displaystyle \sum_{\ell=0}^{\infty}}\dfrac{e^{-1}1^{\ell-1}}{\ell!}\dfrac{\left(-1\right)^{\ell}\left[\left(\ell+1\right)-1\right]}{\left(\ell+1\right)^{2}}
\]
\[
=-{\displaystyle \sum_{\ell=0}^{\infty}}\dfrac{e^{-1}\left(-1\right)^{\ell}}{\left(\ell+1\right)!}+{\displaystyle \sum_{\ell=0}^{\infty}}\dfrac{e^{-1}}{\ell!}\dfrac{\left(-1\right)^{\ell}}{\left(\ell+1\right)^{2}}
\]
\[
=e^{-1}\left(e^{-1}-1\right)+\mu_{-1,1}
\]
\[
=\dfrac{\delta+e\gamma-e+1}{e^{2}}.
\]
This process then can be repeated to solve for $\mu_{-1,3},$ $\mu_{-1,4}$,
etc., yielding the general expression
\begin{equation}
\mu_{-1,k}=\dfrac{\left(k-1\right)!\left(\delta+e\gamma\right)-eA_{k-1}+B_{k-1}}{e^{2}},
\end{equation}
for $k\in\left\{ 2,3,\ldots\right\} $, where: $\left\{ A_{m}\right\} =0,1,3,11,50,274,\ldots$\footnote{Sequence A000254 in the \emph{Online Encyclopedia of Integer Sequences}
(https://oeis.org/). } denotes the sequence of unsigned Stirling numbers of the first kind,
\[
A_{m}=m!{\displaystyle \sum_{i=1}^{m}\dfrac{1}{i}},\quad m\in\left\{ 0,1,\ldots\right\} ;
\]
and $\left\{ B_{m}\right\} =0,1,4,17,84,485,\ldots$\footnote{Sequence A093344 in the \emph{Online Encyclopedia of Integer Sequences}
(https://oeis.org/).} is an unnamed sequence defined by
\begin{equation}
B_{m}=m!{\displaystyle \sum_{i=1}^{m}}\left(\dfrac{1}{i}{\displaystyle \sum_{j=0}^{i-1}}\dfrac{1}{j!}\right),\quad m\in\left\{ 0,1,\ldots\right\} .
\end{equation}

For the limit in (1), it is known that ${\textstyle \sum_{i=1}^{k}}\tfrac{1}{i}=\ln\left(k\right)+\gamma+O\left(\tfrac{1}{k}\right)$
as $k\rightarrow\infty$. This implies
\[
A_{k}=k!\left[\ln\left(k\right)+\gamma+O\left(\dfrac{1}{k}\right)\right]
\]
\begin{equation}
=k!\left(\ln\left(k\right)+\gamma\right)+O\left(\left(k-1\right)!\right),
\end{equation}
from which the asymptotic behavior of $\mu_{-1,k}$, $B_{k}$, and
$eA_{k}-B_{k}$ can be derived.\medskip{}

\noindent \textbf{Proposition 1:} As $k\rightarrow\infty$:

\noindent 
\[
\textrm{(i) }\mu_{-1,k}=O\left(\dfrac{1}{k^{2}}\right);
\]
\[
\textrm{(ii) }B_{k}=ek!\left(\ln\left(k\right)-\dfrac{\delta}{e}\right)+O\left(\left(k-1\right)!\right);
\]

\noindent and

\noindent 
\[
\textrm{(iii) }\underset{k\rightarrow\infty}{\lim}\dfrac{\left(eA_{k}-B_{k}\right)}{e^{2}k!}=\dfrac{\delta+e\gamma}{e^{2}}
\]
\[
=\mu_{-1,1}.
\]
\medskip{}
\textbf{Proof:} (i) From the definition in (6), one can see that
\[
\mu_{-1,k}={\displaystyle \sum_{i=0}^{\infty}}\dfrac{e^{-1}\left(-1\right)^{i}}{i!\left(i+k\right)^{2}}
\]
\[
=\dfrac{1}{k^{2}}{\displaystyle \sum_{i=0}^{\infty}}\dfrac{e^{-1}\left(-1\right)^{i}}{i!}\left(1-\dfrac{i}{k+i}\right)^{2}
\]
\[
=\dfrac{1}{k^{2}}{\displaystyle \sum_{i=0}^{\infty}}\dfrac{e^{-1}\left(-1\right)^{i}}{i!}\left(1+O\left(\dfrac{1}{k}\right)\right)
\]
\[
=\dfrac{e^{-2}}{k^{2}}\left(1+O\left(\dfrac{1}{k}\right)\right)
\]
\[
=O\left(\dfrac{1}{k^{2}}\right).
\]

\noindent (ii) Substituting from (10) into (8) gives
\[
\mu_{-1,k}=\dfrac{\left(k-1\right)!\left(\delta+e\gamma\right)-e\left(k-1\right)!\left[\ln\left(k-1\right)+\gamma+O\left(\dfrac{1}{k-1}\right)\right]+B_{k-1}}{e^{2}}
\]
\[
\Longrightarrow O\left(\dfrac{1}{k^{2}}\right)=\left(k-1\right)!\left(\delta+e\gamma\right)-e\left(k-1\right)!\left[\ln\left(k-1\right)+\gamma+O\left(\dfrac{1}{k-1}\right)\right]+B_{k-1}
\]
\[
\Longrightarrow B_{k-1}=-\left(k-1\right)!\left(\delta+e\gamma\right)+\left(k-1\right)!\left[e\ln\left(k-1\right)+e\gamma+eO\left(\dfrac{1}{k-1}\right)\right]+O\left(\dfrac{1}{k^{2}}\right)
\]
\[
=-\left(k-1\right)!\delta+\left(k-1\right)!\left[e\ln\left(k-1\right)+O\left(\dfrac{1}{k-1}\right)\right]+O\left(\dfrac{1}{k^{2}}\right)
\]
\[
=-\left(k-1\right)!\delta+e\left(k-1\right)!\ln\left(k-1\right)+O\left(\left(k-2\right)!\right),
\]
which implies the desired result.

\noindent (iii) From (8), it follows that
\[
\underset{k\rightarrow\infty}{\lim}\dfrac{\left(eA_{k}-B_{k}\right)}{e^{2}k!}=\underset{k\rightarrow\infty}{\lim}\left[-\dfrac{\mu_{-1,k+1}}{k!}+\dfrac{k!\left(\delta+e\gamma\right)}{e^{2}k!}\right]
\]
\[
=0+\dfrac{\delta+e\gamma}{e^{2}}.\;\blacksquare
\]

\subsection{General $\boldsymbol{k}$ for $\boldsymbol{u=1}$}

\noindent Inserting $u=1$ into the power series (3) yields
\[
\textrm{Ei}\left(1\right)=\gamma+{\displaystyle \sum_{i=1}^{\infty}\dfrac{1^{i}}{i\cdot i!}}
\]
\[
\Longleftrightarrow\delta^{*}=-e\left(\gamma+{\displaystyle \sum_{i=1}^{\infty}\dfrac{1^{i+1}}{i\cdot i!}}\right)
\]
\begin{equation}
\Longleftrightarrow{\displaystyle \sum_{i=1}^{\infty}\dfrac{1^{i+1}}{i\cdot i!}}=-\left(\dfrac{\delta^{*}}{e}+\gamma\right),
\end{equation}
where we define $\delta^{*}\equiv-e\textrm{Ei}\left(1\right)=-5.151464\ldots$
as an analogue of the Euler-Gompertz constant based on the non-alternating
counterpart of the infinite series in (4). This implies
\begin{equation}
\mu_{1,1}=-\dfrac{\left(\delta^{*}+e\gamma\right)}{e^{2}},
\end{equation}
which generalizes to
\[
\mu_{1,k}=\left(-1\right)^{k}\left[\dfrac{\left(k-1\right)!\left(\delta^{*}+e\gamma\right)-eA_{k-1}+e^{2}B_{k-1}^{*}}{e^{2}}\right]
\]
for $k\in\left\{ 2,3,\ldots\right\} $ by iterative calculations similar
to those used to obtain (8). The integer sequence $\left\{ B_{m}^{*}\right\} =0,1,2,7,30,159,\ldots$\footnote{\noindent Sequence A381681 in the \emph{Online Encyclopedia of Integer
Sequences }(https://oeis.org/).} is defined by
\begin{equation}
B_{m}^{*}=m!{\displaystyle \sum_{i=1}^{m}}\left(\dfrac{1}{i}{\displaystyle \sum_{j=0}^{i-1}}\dfrac{\left(-1\right)^{j}}{j!}\right),\quad m\in\left\{ 0,1,\ldots\right\} ,
\end{equation}
which, apart from its alternating component, is identical to (9).\footnote{\noindent The clear similarities of (13) to (9), (12) to (7), and
(11) and (4) suggest that $\delta^{*}$ possesses an ``alternating
analogue'' relationship with $\delta$ comparable to that of $\ln\left(\tfrac{4}{\pi}\right)=0.241564\ldots$
with $\gamma$ (see Sondow, 2005).}

In this case, (10) allows us to derive the following results describing
the asymptotic behavior of $\mu_{1,k}$, $B_{k}^{*}$, and $eA_{k}-e^{2}B_{k}^{*}$.\medskip{}

\noindent \textbf{Proposition 2:} As $k\rightarrow\infty$:

\noindent 
\[
\textrm{(i) }\mu_{1,k}=O\left(\dfrac{1}{k^{2}}\right);
\]
\[
\textrm{(ii) }B_{k}^{*}=\dfrac{k!}{e}\left(\ln\left(k\right)-\dfrac{\delta^{*}}{e}\right)+O\left(\left(k-1\right)!\right);
\]

\noindent and

\noindent 
\[
\textrm{(iii) }\underset{k\rightarrow\infty}{\lim}\dfrac{\left(eA_{k}-e^{2}B_{k}^{*}\right)}{e^{2}k!}=\dfrac{\delta^{*}+e\gamma}{e^{2}}
\]
\[
=-\mu_{1,1}.
\]
\medskip{}
\textbf{Proof:} Parts (i), (ii), and (iii) follow from arguments analogous
to those used to prove the corresponding parts of Proposition 1. $\blacksquare$
\noindent \begin{flushright}
\medskip{}
\par\end{flushright}

\section{A Sequence of Probabilities}

\noindent We now turn to a second conjunction of the Euler constants,
which requires more contextual development than Hardy's equation.
This setting also involves Poisson random variables with $\lambda=1$,
here generated by the continuous-time process $N\left(t\right)\mid\lambda\sim\textrm{Poisson}\left(\lambda t=1\right)$
for $t\in\mathbb{R}^{+}$. For this Poisson process, the random inter-arrival
times, $X_{1}=\min\left(t:N\left(t\right)\geq1\right),X_{2}=\min\left(t:N\left(t\right)\geq2\right)-X_{1},\ldots$
are independent and identically distributed (IID) $\textrm{Exponential}\left(\lambda=1\right)$,
where the parameter $\lambda\in\mathbb{R}^{+}$ denotes the inverse
of the Exponential mean (i.e., $X_{n}\mid\lambda\sim F_{T}\left(t\right)=1-e^{-\lambda t},\:t\in\mathbb{R}^{+}$).
We consider the sequence of random variables, $Y_{n}$, generated
recursively by setting: (1) $Y_{1}\mid\lambda_{1}\equiv X_{1}\mid\lambda_{1}\sim\textrm{Exponential}\left(\lambda_{1}=1\right)$;
and (2) $Y_{n}\mid\lambda_{n}\sim\textrm{Exponential}\left(\lambda_{n}=Y_{n-1}\right)$
for $n\in\left\{ 2,3,\ldots\right\} $.

Unconditionally, the $Y_{n}$ may be expressed as:
\[
Y_{1}\sim F_{Y_{1}}\left(y\right)=1-e^{-y};
\]
\[
Y_{2}\sim F_{Y_{2}}\left(y\right)=E_{Y_{1}}\left[1-e^{-Y_{1}y}\right];
\]
\[
Y_{3}\sim F_{Y_{3}}\left(y\right)=E_{Y_{1}}\left[E_{Y_{2}\mid Y_{1}}\left[1-e^{-Y_{2}y}\right]\right];
\]
\[
Y_{4}\sim F_{Y_{4}}\left(y\right)=E_{Y_{1}}\left[E_{Y_{2}\mid Y_{1}}\left[E_{Y_{3}\mid Y_{2},Y_{1}}\left[1-e^{-Y_{3}y}\right]\right]\right];
\]
\begin{equation}
\textrm{etc}.
\end{equation}
Such models, usually with fairly small values of $n$, can serve as
simple illustrations of parameter uncertainty. In actuarial finance,
for example, they may be used to show the impact of risk heterogeneity
on an insurance loss, $Y_{2}\mid\lambda_{2}\sim\textrm{Exponential}\left(\lambda_{2}\right)$,\footnote{Although the choice of an Exponential model rarely is motivated \foreignlanguage{american}{by
a specific physical process, one readily can imagine that, at the
moment a loss event begins, a continuous impulse of destruction is
applied to an exposed person, property, or other item of value until
some randomly occurring Poisson event disrupts the impulse.} } generated by a single member of a continuum of exposure units with
distinctly different mean losses, $\tfrac{1}{\lambda_{2}}$. Assuming
the particular exposure unit responsible for $Y_{2}$ is selected
randomly \textendash{} and in such a way that $\lambda_{2}\mid\lambda\sim\textrm{Exponential}\left(\lambda=1\right)$
\textendash{} then yields $Y_{2}\sim F_{Y_{2}}\left(y\right)=1-\tfrac{1}{y+1},\:y\in\mathbb{R}^{+}$,
the cumulative distribution function (CDF) of the $\textrm{Pareto 2}\left(\alpha=1,\theta=1\right)$
distribution.

Given the nested conditional expected values presented in (14), one
can rewrite the unconditional CDFs, $F_{Y_{n}}\left(y\right)$, as
follows:
\[
F_{Y_{1}}\left(y\right)=1-e^{-y};
\]
\[
F_{Y_{2}}\left(y\right)={\displaystyle \int_{0}^{\infty}}\left(1-e^{-y_{1}y}\right)e^{-y_{1}}dy_{1}=1-\dfrac{1}{y+1};
\]
\[
F_{Y_{3}}\left(y\right)={\displaystyle \int_{0}^{\infty}}\left[{\displaystyle \int_{0}^{\infty}}\left(1-e^{-y_{2}y}\right)y_{1}e^{-y_{1}y_{2}}dy_{2}\right]e^{-y_{1}}dy_{1}=-ye^{y}\textrm{Ei}\left(-y\right);
\]
\[
F_{Y_{4}}\left(y\right)={\displaystyle \int_{0}^{\infty}}\left[{\displaystyle \int_{0}^{\infty}}\left[{\displaystyle \int_{0}^{\infty}}\left(1-e^{-y_{3}y}\right)y_{2}e^{-y_{2}y_{3}}dy_{3}\right]y_{1}e^{-y_{1}y_{2}}dy_{2}\right]e^{-y_{1}}dy_{1}=\dfrac{y\left(y-\ln\left(y\right)-1\right)}{\left(y-1\right)^{2}};
\]
\[
\textrm{etc}.
\]
Unfortunately, there is no clear sequence of analytic forms for $F_{Y_{n}}\left(y\right)$
as $n$ increases, and the expressions for odd $n$ rapidly become
quite complex (as can be seen from $F_{Y_{3}}\left(y\right)$, which
involves the exponential-integral function). Nevertheless, the recursive
derivation of $Y_{n}$ permits a convenient reformulation. Specifically,
we see that
\[
F_{Y_{2}}\left(y\right)=E_{Y_{1}}\left[1-e^{-Y_{1}y}\right]
\]
\[
=\Pr\left\{ X_{2}\leq Y_{1}y\right\} 
\]
\[
=\Pr\left\{ \dfrac{X_{2}}{Y_{1}}\leq y\right\} ,
\]
from which it follows that $Y_{2}\equiv\tfrac{X_{2}}{Y_{1}}\equiv\tfrac{X_{2}}{X_{1}}$.
Similar arguments yield $Y_{3}\equiv\tfrac{X_{3}}{Y_{2}}\equiv\tfrac{X_{3}X_{1}}{X_{2}}$,
$Y_{4}\equiv\tfrac{X_{4}}{Y_{3}}\equiv\tfrac{X_{4}X_{2}}{X_{3}X_{1}},$
and so on, implying the general identities
\begin{equation}
Y_{n}\equiv\begin{cases}
{\displaystyle \prod_{i=1,3,\ldots,n}}X_{i}\:\:/{\displaystyle \prod_{i=2,4,\ldots,n-1}}X_{i} & \textrm{for }n\in\left\{ 1,3,\ldots\right\} \\
{\displaystyle \prod_{i=2,4,\ldots,n}}X_{i}\:\:/{\displaystyle \prod_{i=1,3,\ldots,n-1}}X_{i} & \textrm{for }n\in\left\{ 2,4,\ldots\right\} 
\end{cases}.
\end{equation}

The expressions for $Y_{n}$ in (15) reveal two important aspects
of this sequence. First, the random variable $Y_{n}$ possesses the
same probability distribution as its inverse ($Y_{n}^{-1}$) for all
even $n$, with both distributions ``symmetric'' about $y=1$ in
the sense that $\Pi_{n}=F_{Y_{n}}\left(1\right)=F_{Y_{n}^{-1}}\left(1\right)=\tfrac{1}{2}$.
Second, $Y_{n}$ is somewhat ``top-heavy'' for odd $n$, with $\Pi_{n}=F_{Y_{n}}\left(1\right)>\tfrac{1}{2}$.

Computed values of $\Pi_{n}$ provided by Table 1 show that the impact
of the additional Exponential random variable in the numerator of
$Y_{n}$ decreases over $n$, with $\Pi_{n}\rightarrow\tfrac{1}{2}$
as $n\rightarrow\infty$. Moreover, the table reveals an additional,
somewhat intriguing, property: each of the probabilities $\Pi_{1}$,
$\Pi_{3}$, and $\Pi_{5}$ is a simple function of one of the three
Euler constants, $e$, $\delta$, and $\gamma$, respectively.\pagebreak{}
\noindent \begin{center}
Table 1. Values of $\Pi_{n}=F_{Y_{n}}\left(1\right),\:n\in\left\{ 1,2,\ldots,10\right\} $
\par\end{center}

\noindent \begin{center}
\begin{tabular}{|c|c|}
\hline 
$n$ & $\Pi_{n}$\tabularnewline
\hline 
\hline 
$1$ & $0.632120\ldots=1-e^{-1}$\tabularnewline
\hline 
$2$ & $0.5$\tabularnewline
\hline 
$3$ & $0.596347\ldots=\delta$\tabularnewline
\hline 
$4$ & $0.5$\tabularnewline
\hline 
$5$ & $0.577215\ldots=\gamma$\tabularnewline
\hline 
$6$ & $0.5$\tabularnewline
\hline 
$7$ & $0.566094\ldots$\tabularnewline
\hline 
$8$ & $0.5$\tabularnewline
\hline 
$9$ & $0.558672\ldots$\tabularnewline
\hline 
$10$ & $0.5$\tabularnewline
\hline 
\end{tabular}
\par\end{center}

\medskip{}

\subsection{Log-Scale Analysis}

\noindent As noted above, the random variable $Y_{2}$ (which can
be interpreted as the ratio of two independent $\textrm{Exponential}\left(\lambda=1\right)$
random variables), possesses a $\textrm{Pareto 2}\left(\alpha=1,\theta=1\right)$
distribution, and therefore an infinite mean. As a result, all $Y_{n}$
for $n>2$ are similarly heavy-tailed, which hinders closer examination
of their distributional properties through moment calculations and
the central limit theorem (CLT). For that reason, we now transform
$Y_{n}$ to the log scale, working with $Z_{n}\equiv\ln\left(Y_{n}\right)$
for $n\in\left\{ 1,2,\ldots\right\} $.

\subsubsection{Asymptotic Distributions}

\noindent In addition to sidestepping the problem of heavy tails,
the log transformation enables us to work with sums, rather than products,
of independent random variables. In particular, we can write
\begin{equation}
Z_{n}\equiv\begin{cases}
{\displaystyle \sum_{i=1,3,\ldots,n}}\ln\left(X_{i}\right)-{\displaystyle \sum_{i=2,4,\ldots,n-1}}\ln\left(X_{i}\right) & \textrm{for }n\in\left\{ 1,3,\ldots\right\} \\
{\displaystyle \sum_{i=2,4,\ldots,n}}\ln\left(X_{i}\right)-{\displaystyle \sum_{i=1,3,\ldots,n-1}}\ln\left(X_{i}\right) & \textrm{for }n\in\left\{ 2,4,\ldots\right\} 
\end{cases},
\end{equation}
where the $-\ln\left(X_{i}\right)\sim\textrm{IID Gumbel}\left(m=0,s=1\right)$
(with mean $\gamma$ and variance $\tfrac{\pi^{2}}{6}$; i.e., $-\ln\left(X_{i}\right)\equiv U_{i}\sim F_{U}\left(u\right)=\exp\left(-e^{-u}\right),\:u\in\mathbb{R}$).\footnote{For $n\in\left\{ 2,4,\ldots\right\} $, this implies $Z_{n}$ is the
$\tfrac{n}{2}$-fold convolution of IID $\textrm{Logistic}\left(m=0,s=1\right)$
random variables; i.e., $Z_{n}={\textstyle \sum_{j=1}^{n/2}}V_{j}$
for $V_{j}\sim\textrm{IID }F_{V}\left(v\right)=\tfrac{1}{1+\exp\left(-v\right)},\:v\in\mathbb{R}$.
Hereafter, we will say $Z_{n}\sim\tfrac{n}{2}\textrm{-Fold Logistic}\left(0,1\right)$
and $Y_{n}\sim\textrm{Log }\tfrac{n}{2}\textrm{-Fold Logistic}\left(0,1\right)$
for $n\in\left\{ 2,4,\ldots\right\} $.} A straightforward application of the CLT then yields the following
result.\medskip{}

\noindent \textbf{Proposition 3:} For the sequences of random variables,
$Y_{n}$ and $Z_{n}$ (defined in (15) and (16), respectively),
\[
\dfrac{\ln\left(Y_{n}\right)}{\pi\sqrt{\dfrac{n}{6}}}\equiv\dfrac{Z_{n}}{\pi\sqrt{\dfrac{n}{6}}}\stackrel{\mathcal{D}}{\longrightarrow}\begin{cases}
\textrm{Normal}\left(-\dfrac{\gamma}{\pi\sqrt{\dfrac{n}{6}}},\:1\right) & \textrm{for }n\in\left\{ 1,3,\ldots\right\} \\
\textrm{Normal}\left(0,\:1\right) & \textrm{for }n\in\left\{ 2,4,\ldots\right\} 
\end{cases}
\]
as $n\rightarrow\infty$.\medskip{}

\noindent \textbf{Proof:} For $n\in\left\{ 1,3,\ldots\right\} $,
it is clear from (16) that
\[
Z_{n}\equiv-{\displaystyle \sum_{i=1,3,\ldots,n}U_{i}}+{\displaystyle \sum_{i=2,4,\ldots,n-1}U_{i}},
\]
where the $U_{i}$ are IID $\textrm{Gumbel}\left(m=0,s=1\right)$
random variables. Letting $Z_{n}^{\left(1\right)}=-{\displaystyle {\textstyle \sum_{i=1,3,\ldots,n}}U_{i}}$
and $Z_{n}^{\left(2\right)}={\displaystyle {\textstyle \sum_{i=2,4,\ldots,n-1}}U_{i}}$,
it follows from the CLT that $\tfrac{Z_{n}^{\left(1\right)}}{\pi\sqrt{\left(n+1\right)/12}}\stackrel{\mathcal{D}}{\longrightarrow}\textrm{Normal}\left(-\tfrac{\gamma}{\pi\sqrt{\left(n+1\right)/12}},1\right)$
and $\tfrac{Z_{n}^{\left(2\right)}}{\pi\sqrt{\left(n-1\right)/12}}\stackrel{\mathcal{D}}{\longrightarrow}\textrm{Normal}\left(0,1\right)$.
Taking the linear combination
\[
Z_{n}=\pi\sqrt{\dfrac{n+1}{12}}Z_{n}^{\left(1\right)}+\pi\sqrt{\dfrac{n-1}{12}}Z_{n}^{\left(2\right)}
\]
then gives the desired result.

For $n\in\left\{ 2,4,\ldots\right\} $,
\[
Z_{n}\equiv-{\displaystyle \sum_{i=2,4,\ldots,n}^{n/2}U_{i}}+{\displaystyle \sum_{i=1,3,\ldots,n-1}^{n/2}U_{i}},
\]
and the argument proceeds in the same way. $\blacksquare$\medskip{}

The above proposition thus reveals that $Y_{n}$ is asymptotically
Lognormal with increasing accumulations of the total probability split
equally between two regions: the distant right tail, and a small neighborhood
of 0.

\subsubsection{Characteristic Functions and CDFs}

\noindent Another benefit of working with sums of IID random variables
is that expressions for the characteristic function are more likely
to be tractable. In the case at hand, we obtain the following proposition.\medskip{}

\noindent \textbf{Proposition 4:} For the sequence of random variables,
$Z_{n}$, the corresponding characteristic functions are given by:
\begin{equation}
\textrm{(i) }\varphi_{Z_{n}}\left(\omega\right)=\begin{cases}
\left[\Gamma\left(1+i\omega\right)\right]^{\left(n+1\right)/2}\left[\Gamma\left(1-i\omega\right)\right]^{\left(n-1\right)/2} & \textrm{for }n\in\left\{ 1,3,\ldots\right\} \\
\left[\Gamma\left(1+i\omega\right)\right]^{n/2}\left[\Gamma\left(1-i\omega\right)\right]^{n/2} & \textrm{for }n\in\left\{ 2,4,\ldots\right\} 
\end{cases};
\end{equation}
or equivalently,
\begin{equation}
\textrm{(ii) }\varphi_{Z_{n}}\left(\omega\right)=\begin{cases}
\left(\dfrac{\pi\omega}{\sinh\left(\pi\omega\right)}\right)^{\left(n-1\right)/2}{\displaystyle \int_{0}^{\infty}e^{-t}\left[\cos\left(\ln\left(t\right)\omega\right)+i\sin\left(\ln\left(t\right)\omega\right)\right]dt} & \textrm{for }n\in\left\{ 1,3,\ldots\right\} \\
\left(\dfrac{\pi\omega}{\sinh\left(\pi\omega\right)}\right)^{n/2} & \textrm{for }n\in\left\{ 2,4,\ldots\right\} 
\end{cases}.
\end{equation}
\medskip{}
\textbf{Proof:} From (16), it is easy to see that
\begin{equation}
\varphi_{Z_{n}}\left(\omega\right)=\begin{cases}
\left[\varphi_{-\ln\left(X\right)}\left(-\omega\right)\right]^{\left(n+1\right)/2}\left[\varphi_{-\ln\left(X\right)}\left(\omega\right)\right]^{\left(n-1\right)/2} & \textrm{for }n\in\left\{ 1,3,\ldots\right\} \\
\left[\varphi_{-\ln\left(X\right)}\left(-\omega\right)\right]^{n/2}\left[\varphi_{-\ln\left(X\right)}\left(\omega\right)\right]^{n/2} & \textrm{for }n\in\left\{ 2,4,\ldots\right\} 
\end{cases},
\end{equation}
where $-\ln\left(X\right)\sim\textrm{Gumbel}\left(m=0,s=1\right)$.
Since $\varphi_{-\ln\left(X\right)}\left(-\omega\right)=E_{X}\left[e^{i\omega\ln\left(X\right)}\right]=E_{X}\left[X^{i\omega}\right]=\Gamma\left(1+i\omega\right)$,
we may obtain (17) from direct substitution into (19), and then (18)
from Euler's reflection formula and the identity $\Gamma\left(1+i\omega\right)\equiv\int_{0}^{\infty}e^{-t}t^{i\omega}dt$.
$\blacksquare$\medskip{}

Rather fortuitously, the characteristic functions of Proposition 4
are relatively easy to invert, as shown in the next result.\medskip{}

\noindent \textbf{Proposition 5:} For the sequence of random variables,
$Z_{n}$, the corresponding CDFs are given by:
\begin{equation}
\textrm{(i) }F_{Z_{n}}\left(z\right)=\begin{cases}
{\displaystyle G_{1,2}^{1,1}}\left(e^{z}\left|\begin{array}{c}
1\\
1,0
\end{array}\right.\right) & \textrm{for }n=1\\
{\displaystyle G_{\left(n-1\right)/2,\left(n+1\right)/2}^{\left(n+1\right)/2,\left(n-1\right)/2}}\left(e^{z}\left|\begin{array}{c}
1,0,\ldots,0\\
1,1,\ldots,1,1
\end{array}\right.\right) & \textrm{for }n\in\left\{ 3,5,\ldots\right\} \\
{\displaystyle G_{n/2,n/2}^{n/2,n/2}}\left(e^{z}\left|\begin{array}{c}
1,0,\ldots,0\\
1,1,\ldots,1
\end{array}\right.\right) & \textrm{for }n\in\left\{ 2,4,\ldots\right\} 
\end{cases},
\end{equation}
where $G$ denotes the Meijer G-function (as parameterized by Mathematica);
and equivalently,
\begin{equation}
\textrm{(ii) }F_{Z_{n}}\left(z\right)=\begin{cases}
\dfrac{1}{2}+{\displaystyle \int_{0}^{\infty}}{\displaystyle \int_{0}^{\infty}}e^{-t}\dfrac{\sin\left(\left[z-\ln\left(t\right)\right]\omega\right)}{\pi\omega}{\displaystyle \left(\dfrac{\pi\omega}{\sinh\left(\pi\omega\right)}\right)^{\left(n-1\right)/2}}dt{\displaystyle d\omega} & \textrm{for }n\in\left\{ 1,3,\ldots\right\} \\
\dfrac{1}{2}+{\displaystyle \int_{0}^{\infty}\dfrac{\sin\left(z\omega\right)}{\pi\omega}\left(\dfrac{\pi\omega}{\sinh\left(\pi\omega\right)}\right)^{n/2}{\displaystyle d\omega}} & \textrm{for }n\in\left\{ 2,4,\ldots\right\} 
\end{cases};
\end{equation}
for $z\in\mathbb{R}$.\medskip{}

\noindent \textbf{Proof:} See the Appendix.\medskip{}

Since $F_{Y_{n}}\left(1\right)=F_{Z_{n}}\left(0\right)$ for all $n$,
we may obtain expressions for the probabilities $\Pi_{n}$, for $n\in\left\{ 1,3,\ldots\right\} $,
by setting $z=0$ in (20) and (21). This yields
\begin{equation}
\Pi_{n}=\begin{cases}
{\displaystyle G_{1,2}^{1,1}}\left(1\left|\begin{array}{c}
1\\
1,0
\end{array}\right.\right) & \textrm{for }n=1\\
{\displaystyle G_{\left(n-1\right)/2,\left(n+1\right)/2}^{\left(n+1\right)/2,\left(n-1\right)/2}}\left(1\left|\begin{array}{c}
1,0,\ldots,0\\
1,1,\ldots,1,1
\end{array}\right.\right) & \textrm{for }n\in\left\{ 3,5,\ldots\right\} 
\end{cases}
\end{equation}
and
\begin{equation}
\Pi_{n}=\dfrac{1}{2}-\int_{0}^{\infty}\int_{0}^{\infty}e^{-t}\dfrac{\sin\left(\ln\left(t\right)\omega\right)}{\pi\omega}{\displaystyle \left(\dfrac{\pi\omega}{\sinh\left(\pi\omega\right)}\right)^{\left(n-1\right)/2}}dt{\displaystyle d\omega},
\end{equation}
respectively.

Although the above mathematical forms provide only limited analytic
tractability, they do afford expressions for $\Pi_{n}$ as explicit
functions of $n$ that may be helpful in certain contexts. For example,
to compute $\Pi_{n}$ to a high degree of accuracy, one can take advantage
of efficient numerical algorithms for the Meijer G-Function in (22).
As an illustration, we present the results of such calculations to
50 decimal places (using Mathematica) in Table 2 below. Moreover,
the integral in (23), although not as computationally convenient as
(22), offers the advantage of embedding $\Pi_{n}$ \textendash{} and
in particular, $\delta$ and $\gamma$ \textendash{} into a continuous
function of the parameter $n$. This may be useful in analytic studies
of these quantities.\medskip{}
\noindent \begin{center}
Table 2. Values of $\Pi_{n},\:n\in\left\{ 1,3,\ldots,15\right\} $
Computed by Mathematica
\par\end{center}

\noindent \begin{center}
\begin{tabular}{|c|c|}
\hline 
$n$ & $\Pi_{n}$\tabularnewline
\hline 
\hline 
$1$ & $0.63212055882855767840447622983853913255418886896823\ldots$\tabularnewline
\hline 
$3$ & $0.59634736232319407434107849936927937607417786015254\ldots$\tabularnewline
\hline 
$5$ & $0.57721566490153286060651209008240243104215933593992\ldots$\tabularnewline
\hline 
$7$ & $0.56609435541264796901908591583288674247188413864361\ldots$\tabularnewline
\hline 
$9$ & $0.55867279019459907350395199904241483559079945290197\ldots$\tabularnewline
\hline 
$11$ & $0.55328267668997479292771549900978480433113900145274\ldots$\tabularnewline
\hline 
$13$ & $0.54914332831601761785791255217329440945377229653973\ldots$\tabularnewline
\hline 
$15$ & $0.54583694813712457806697974073677754090116258346196\ldots$\tabularnewline
\hline 
\end{tabular}\medskip{}
\par\end{center}

\subsection{Closed-Form Expressions}

\noindent Further analysis facilitates the construction of closed-form
expressions for the $\Pi_{n}$, as we now show.

\subsubsection{The $\boldsymbol{\tfrac{n}{2}-\textrm{Fold Logistic}\left(0,1\right)}$
Survival Function}

\noindent From (15), it is known that $Y_{n}\equiv\tfrac{X_{n}}{Y_{n-1}}$,
implying
\[
F_{Y_{n}}\left(1\right)=\Pr\left\{ \dfrac{X_{n}}{Y_{n-1}}\leq1\right\} 
\]
\[
=\Pr\left\{ X_{n}\leq Y_{n-1}\right\} 
\]
\[
=E_{X_{n}}\left[\Pr\left\{ \left.X_{n}\leq Y_{n-1}\right|X_{n}\right\} \right].
\]
For $n\in\left\{ 3,5,\ldots\right\} $, we then can write
\[
\Pi_{n}{\displaystyle =}\int_{0}^{\infty}\left[1-F_{Y_{n-1}}\left(x\right)\right]e^{-x}dx
\]
\[
{\displaystyle =}\int_{0}^{\infty}S_{Y_{n-1}}\left(x\right)e^{-x}dx
\]
\[
{\displaystyle =}\int_{0}^{\infty}S_{Z_{n-1}}\left(\ln\left(x\right)\right)e^{-x}dx,
\]
where $S_{Y_{n-1}}\left(\cdot\right)$ and $S_{Z_{n-1}}\left(\cdot\right)$
denote the respective survival functions of the $\textrm{Log }\tfrac{n-1}{2}\textrm{-Fold Logistic}\left(0,1\right)$
and $\tfrac{n-1}{2}\textrm{-Fold Logistic}\left(0,1\right)$ distributions,
for $n-1\in\left\{ 2,4,\ldots\right\} $.\footnote{See Footnote 6.}

Although analytic expressions for the CDF and PDF of $Z_{n-1}$ have
appeared in the literature (see George and Mudholkar, 1983 and Ojo,
1986, respectively), they are not well-suited for the purpose at hand.
Instead, we will work with the following compact form of the survival
function (evaluated at $\ln\left(y\right)$):
\begin{equation}
S_{Z_{n-1}}\left(\ln\left(y\right)\right)=\dfrac{{\displaystyle \sum_{i=1}^{\left(n-3\right)/2}}\:{\displaystyle \sum_{j=0}^{\left(n-3\right)/2}}C_{i,j}^{\left(n-1\right)}y^{i}\left(\ln\left(y\right)\right)^{j}+1}{\left[y-\left(-1\right)^{\left(n-1\right)/2}\right]^{\left(n-1\right)/2}},\:n-1\in\left\{ 2,4,\ldots\right\} ,
\end{equation}
where the $C_{i,j}^{\left(n-1\right)}$ are real-valued coefficients.

Rather than deriving (24) directly from the second line of (21) (which
is quite tedious), we demonstrate its validity by showing that it
satisfies equivalent conditions. From (15), we know that $Y_{n-1}\equiv\tfrac{X_{n-1}}{X_{n-2}}Y_{n-3}$
for $n-1\in\left\{ 4,6,\ldots\right\} $, yielding
\[
S_{Y_{n-1}}\left(y\right)=\Pr\left\{ \dfrac{X_{n-1}}{X_{n-2}}Y_{n-3}>y\right\} 
\]
\[
=\Pr\left\{ Y_{n-3}>\dfrac{X_{n-2}}{X_{n-1}}y\right\} 
\]
\[
=E_{X_{n-1},X_{n-2}}\left[\Pr\left\{ \left.Y_{n-3}>\dfrac{X_{n-2}}{X_{n-1}}y\right|X_{n-1},X_{n-2}\right\} \right],
\]
\[
{\displaystyle =}\int_{0}^{\infty}\dfrac{S_{Y_{n-3}}\left(\tau y\right)}{\left(\tau+1\right)^{2}}d\tau,
\]
where $\tfrac{1}{\left(\tau+1\right)^{2}}$ is the PDF of $\tfrac{X_{n-2}}{X_{n-1}}\sim\textrm{Pareto 2}\left(\alpha=1,\theta=1\right)$.
Substituting $\upsilon=\tau y$ then gives
\[
S_{Y_{n-1}}\left(y\right)=y\int_{0}^{\infty}\dfrac{S_{Y_{n-3}}\left(\upsilon\right)}{\left(\upsilon+y\right)^{2}}d\upsilon
\]
\[
=-y\dfrac{\partial}{\partial y}\int_{0}^{\infty}\dfrac{S_{Y_{n-3}}\left(\upsilon\right)}{\upsilon+y}d\upsilon
\]
\begin{equation}
\Longleftrightarrow S_{Z_{n-1}}\left(\ln\left(y\right)\right)=-y\dfrac{\partial}{\partial y}\int_{0}^{\infty}\dfrac{S_{Z_{n-3}}\left(\ln\left(\upsilon\right)\right)}{\upsilon+y}d\upsilon,
\end{equation}
a condition that uniquely defines $S_{Z_{n-1}}\left(\ln\left(y\right)\right)$
for $n-1\in\left\{ 4,6,\ldots\right\} $ subject to
\[
S_{Z_{2}}\left(\ln\left(y\right)\right)=\dfrac{1}{y+1}
\]
(the survival function of $e^{Z_{2}}\sim\textrm{Pareto 2}\left(\alpha=1,\theta=1\right)$).

Since it is fairly easy to work out
\begin{equation}
S_{Z_{4}}\left(\ln\left(y\right)\right)=\dfrac{y\ln\left(y\right)-y+1}{\left(y-1\right)^{2}}
\end{equation}
from the second line of (21), we will use the more complicated case
of $n-1=6$ to illustrate how one arrives at (24) from (25), noting
that larger values of $n-1$ may be addressed in the same manner with
the assistance of a computer algebra system. Inserting (26) into the
integral in (25) then yields
\[
S_{Z_{6}}\left(\ln\left(y\right)\right)=-y\dfrac{\partial}{\partial y}\int_{0}^{\infty}\dfrac{\upsilon\ln\left(\upsilon\right)-\upsilon+1}{\left(\upsilon-1\right)^{2}\left(\upsilon+y\right)}d\upsilon
\]
\[
=-y\dfrac{\partial}{\partial y}\int_{0}^{\infty}\dfrac{\upsilon\ln\left(\upsilon\right)}{\left(\upsilon-1\right)^{2}\left(\upsilon+y\right)}d\upsilon+y\dfrac{\partial}{\partial y}\int_{0}^{\infty}\dfrac{1}{\left(\upsilon-1\right)\left(\upsilon+y\right)}d\upsilon
\]
\[
=-y\dfrac{\partial}{\partial y}\left[\dfrac{y\left[\pi^{2}+\left(\ln\left(y\right)\right)^{2}\right]}{2\left(y+1\right)^{2}}\right]+y\dfrac{\partial}{\partial y}\left[\dfrac{\ln\left(y\right)}{y+1}\right]
\]
\[
=\dfrac{y\left(y-1\right)\left[\pi^{2}+\left(\ln\left(y\right)\right)^{2}\right]-2y\left(y+1\right)\ln\left(y\right)}{2\left(y+1\right)^{3}}-\dfrac{\left(y\ln\left(y\right)-y-1\right)}{\left(y+1\right)^{2}}
\]
\begin{equation}
=\dfrac{\dfrac{1}{2}\left[\left(\ln\left(y\right)\right)^{2}-4\ln\left(y\right)+\pi^{2}+2\right]y^{2}-\dfrac{1}{2}\left[\left(\ln\left(y\right)\right)^{2}+4\ln\left(y\right)+\pi^{2}-4\right]y+1}{\left(y+1\right)^{3}},
\end{equation}
where both integrals in the second line of (27) are evaluated by the
method of partial fractions, with the former requiring principal-value
cancellation.

\subsubsection{Integration Results}

\noindent The general form of (24) permits the construction of closed-form
expressions for all $\Pi_{n}$ ($n\in\left\{ 3,5,\ldots\right\} $)
via (25). Letting $\zeta\left(\cdot\right)$ denote the Riemann zeta
function and defining
\[
\sigma_{m}={\displaystyle \sum_{i=1}^{\infty}\dfrac{\left(-1\right)^{i+1}}{i^{m}\cdot i!}}
\]
for $m\in\left\{ 1,2,\ldots\right\} $,\footnote{Note that we previously encountered $\sigma_{1}$ in (4) (i.e., Hardy's
equation). } we then have:
\[
\Pi_{3}{\displaystyle =}\int_{0}^{\infty}\left(\dfrac{1}{x+1}\right)e^{-x}dx
\]
\[
=\delta,
\]
which is well known;
\[
\Pi_{5}{\displaystyle =}\int_{0}^{\infty}\left[\dfrac{x\ln\left(x\right)-x+1}{\left(x-1\right)^{2}}\right]e^{-x}dx
\]
\[
=\gamma,
\]
which is possibly new to the literature;\footnote{In particular, the integral does not appear in Choi and Srivastava
(2010) or various online compendiums of expressions for the Euler-Mascheroni
constant. } and
\[
\Pi_{7}{\displaystyle =}\int_{0}^{\infty}S_{Z_{6}}\left(\ln\left(x\right)\right)e^{-x}dx
\]
\[
{\displaystyle =}\dfrac{1}{2}\int_{0}^{\infty}\left[\dfrac{\left(\ln\left(x\right)\right)^{2}x^{2}}{\left(x+1\right)^{3}}\right]e^{-x}dx-2\int_{0}^{\infty}\left[\dfrac{\ln\left(x\right)x^{2}}{\left(x+1\right)^{3}}\right]e^{-x}dx+\dfrac{\left(\pi^{2}+2\right)}{2}\int_{0}^{\infty}\left[\dfrac{x^{2}}{\left(x+1\right)^{3}}\right]e^{-x}dx
\]
\[
-\dfrac{1}{2}\int_{0}^{\infty}\left[\dfrac{\left(\ln\left(x\right)\right)^{2}x}{\left(x+1\right)^{3}}\right]e^{-x}dx-2\int_{0}^{\infty}\left[\dfrac{\ln\left(x\right)x}{\left(x+1\right)^{3}}\right]e^{-x}dx-\dfrac{\left(\pi^{2}-4\right)}{2}\int_{0}^{\infty}\left[\dfrac{x}{\left(x+1\right)^{3}}\right]e^{-x}dx
\]
\[
+\int_{0}^{\infty}\left[\dfrac{1}{\left(x+1\right)^{3}}\right]e^{-x}dx
\]
\begin{equation}
=\dfrac{5}{6}e\left(-\gamma^{3}+3\sigma_{1}\gamma^{2}-21\gamma\zeta\left(2\right)-2\zeta\left(3\right)+21\sigma_{1}\zeta\left(2\right)-6\sigma_{2}\gamma+6\sigma_{3}\right)-\dfrac{3}{2}\left(\gamma^{2}+\dfrac{2}{3}\gamma+7\zeta\left(2\right)\right),
\end{equation}
where all seven integrals in the second line of (28) can be evaluated
using Mellin transforms and the polygamma function.

The expression in (28) is noteworthy because it is dramatically more
complicated than the expressions for $\Pi_{1}$, $\Pi_{3}$, and $\Pi_{5}$.
As in the case of $\Pi_{7}$, each $\Pi_{n},\:n\in\left\{ 9,11,\ldots\right\} $
requires the evaluation of $\tfrac{\left(n-1\right)\left(n-3\right)}{4}+1=\tfrac{n^{2}-4n+7}{4}$
integrals of the form ${\textstyle \int_{0}^{\infty}}\left[\tfrac{x^{i}\left(\ln\left(x\right)\right)^{j}}{\left[x-\left(-1\right)^{\left(n-1\right)/2}\right]^{\left(n-1\right)/2}}\right]e^{-x}dx$
(including ${\textstyle \int_{0}^{\infty}}\left[\tfrac{1}{\left[x-\left(-1\right)^{\left(n-1\right)/2}\right]^{\left(n-1\right)/2}}\right]e^{-x}dx$).
A close inspection of these integrals reveals that the degree of complexity
increases over $n$, with each successive $\Pi_{n}$ expressible as
a linear combination of the constants $e$, $\gamma$, $\zeta\left(2\right)$,
$\zeta\left(3\right)$, \dots , $\zeta\left(\tfrac{n-1}{2}\right)$,
$\sigma_{1}$, $\sigma_{2}$, \dots , $\sigma_{\left(n-1\right)/2}$
(and various products thereof) with rational coefficients.

\subsubsection{Additional Insights; Reprise of Hardy's Equation}

\noindent The constants $\sigma_{m}$ arise occasionally in the literature,
but usually only for certain small values of $m$ or in conjunction
with derivations involving other series. In particular, the following
explicit connection to $U\sim\textrm{Gumbel}\left(0,1\right)$ does
not appear to be widely known:
\[
\sigma_{m}=\dfrac{1}{m!}{\displaystyle \int_{0}^{\infty}}u^{m}\exp\left(-u-e^{-u}\right)du
\]
\[
=\dfrac{1}{m!}{\displaystyle \textrm{E}_{U}\left[\left(U^{+}\right)^{m}\right]},
\]
where $U^{+}=\max\left\{ 0,U\right\} $. This identity, which is readily
demonstrated by substituting $t=e^{-u}$ and applying repeated integration
by parts to the first line, offers an immediate probabilistic interpretation
of Hardy's equation. Since $\textrm{E}\left[U\right]=\gamma$ and
$\Pi_{1}=1-\tfrac{1}{e}=\Pr\left\{ U>0\right\} $, we see that (4)
is equivalent to
\begin{equation}
\textrm{E}_{U}\left[U^{+}\right]=\delta\Pr\left\{ U\leq0\right\} +\textrm{E}_{U}\left[U\right],
\end{equation}
implying
\[
\delta=\dfrac{\textrm{E}_{U}\left[U^{-}\right]}{\Pr\left\{ U\leq0\right\} }
\]
\[
=-\textrm{E}_{U}\left[U\mid U\leq0\right]
\]
\[
=-e{\displaystyle \int_{-\infty}^{0}}u\exp\left(-u-e^{-u}\right)du.
\]

Continuing in this vein, one can define
\[
\delta_{m}=-\textrm{E}_{U}\left[U^{m}\mid U\leq0\right]
\]
\[
=-e{\displaystyle \int_{-\infty}^{0}}u^{m}\exp\left(-u-e^{-u}\right)du
\]
for $m\in\left\{ 1,2,\ldots\right\} $ (where $\delta_{1}\equiv\delta$)
and extend Hardy's equation in a manner very different from that of
Subsections 2.1 and 2.2. Specifically, set
\[
{\displaystyle \textrm{E}_{U}\left[\left(U^{+}\right)^{m}\right]}=\delta_{m}\Pr\left\{ U\leq0\right\} +\textrm{E}_{U}\left[U^{m}\right]
\]
\[
\Longleftrightarrow m!\sigma_{m}=\dfrac{\delta_{m}}{e}+\textrm{E}_{U}\left[U^{m}\right]
\]
(in the forms of (29) and (4), respectively), and note that this gives
\[
2\sigma_{2}=\dfrac{\delta_{2}}{e}+\gamma^{2}+\zeta\left(2\right),
\]
\[
6\sigma_{3}=\dfrac{\delta_{3}}{e}+\gamma^{3}+3\gamma\zeta\left(2\right)+2\zeta\left(3\right),
\]
\[
24\sigma_{4}=\dfrac{\delta_{4}}{e}+\gamma^{4}+6\zeta\left(2\right)\gamma^{2}+8\zeta\left(3\right)\gamma+\dfrac{27}{2}\zeta\left(4\right),
\]
\[
\textrm{etc.}
\]
as analogues to
\[
\sigma_{1}=\dfrac{\delta_{1}}{e}+\gamma
\]
(where expressions for $\textrm{E}_{U}\left[U^{2}\right]$, $\textrm{E}_{U}\left[U^{3}\right]$,
and $\textrm{E}_{U}\left[U^{4}\right]$ are taken from Weisstein,
2025).

Finally, the indicated values of $\sigma_{1}$, $\sigma_{2}$, and
$\sigma_{3}$ may be substituted into (28), yielding
\[
\Pi_{7}=\left(\dfrac{5}{2}\delta_{1}-\dfrac{3}{2}\right)\left(\gamma^{2}+7\zeta\left(2\right)\right)-\left(\dfrac{5}{2}\delta_{2}+1\right)\gamma+\dfrac{5}{6}\delta_{3},
\]
a more concise expression built from 5, rather than 7, distinct constants.
As with $\Pi_{1}$, $\Pi_{3}$, and $\Pi_{5}$, each constant is easily
interpreted in terms of parameters of the $\textrm{Gumbel}\left(0,1\right)$
distribution. 

\section{Conclusion}

\noindent In the present work, we explored two conjunctions of the
three Euler constants ($e$, $\gamma$, and $\delta$), both of which
are associated with the Poisson probability distribution. For the
first case, involving a well-known equation of Hardy (1949), extensions
of this equation were given based on its interpretation in terms of
an inverse second moment of the Poisson random variable. Since the
second conjunction arises from a less familiar context, involving
a sequence of recursively generated Exponential waiting times, we
first lay the groundwork necessary to embed the three Euler constants
into a specific family of probability distributions. Viewing these
constants as salient components of the first three terms in an infinite
sequence of (presumably irrational) numbers, we were able to obtain
new insights into their connections with one another, as well as with
other well-known mathematical quantities. In particular, this analysis
enabled an alternative extension of Hardy's equation.

\section*{Appendix}

\noindent \textbf{Proof of Proposition 5:}

\noindent For part (i), we provide a full derivation for $n\in\left\{ 1,3,\ldots\right\} $,
noting that the simpler case of $n\in\left\{ 2,4,\ldots\right\} $
proceeds analogously. Inserting the first line of (17) into the inversion
formula of Gil-Pelaez (1951) yields
\[
F_{Z_{n}}\left(z\right)=\dfrac{1}{2}+\dfrac{1}{2\pi}{\displaystyle \int_{0}^{\infty}\dfrac{\left[e^{i\omega z}\varphi_{Z_{n}}\left(-\omega\right)-e^{-i\omega z}\varphi_{Z_{n}}\left(\omega\right)\right]}{i\omega}d\omega}\qquad\qquad(\textrm{A1)}
\]
\[
=\dfrac{1}{2}+\dfrac{1}{2\pi}{\displaystyle \int_{0}^{\infty}\dfrac{e^{i\omega z}\left[\Gamma\left(1-i\omega\right)\right]^{\left(n+1\right)/2}\left[\Gamma\left(1+i\omega\right)\right]^{\left(n-1\right)/2}}{i\omega}d\omega}
\]
\[
-\dfrac{1}{2\pi}{\displaystyle \int_{0}^{\infty}\dfrac{e^{-i\omega z}\left[\Gamma\left(1+i\omega\right)\right]^{\left(n+1\right)/2}\left[\Gamma\left(1-i\omega\right)\right]^{\left(n-1\right)/2}}{i\omega}d\omega},
\]
which, after making the substitution $\omega^{\prime}=-\omega$ in
the second integral, becomes
\[
F_{Z_{n}}\left(z\right)=\dfrac{1}{2}+\dfrac{1}{2\pi}{\displaystyle \int_{0}^{\infty}\dfrac{e^{i\omega z}\left[\Gamma\left(1-i\omega\right)\right]^{\left(n+1\right)/2}\left[\Gamma\left(1+i\omega\right)\right]^{\left(n-1\right)/2}}{i\omega}d\omega}
\]
\[
+\dfrac{1}{2\pi}{\displaystyle \int_{-\infty}^{0}\dfrac{e^{i\omega^{\prime}z}\left[\Gamma\left(1-i\omega^{\prime}\right)\right]^{\left(n+1\right)/2}\left[\Gamma\left(1+i\omega^{\prime}\right)\right]^{\left(n-1\right)/2}}{i\omega^{\prime}}d\omega^{\prime}}.
\]
To join the two integrals together on the real line, one can integrate
over an arc contour in the upper half-plane excluding the pole at
the origin, obtaining
\[
F_{Z_{n}}\left(z\right)=\dfrac{1}{2}+\dfrac{1}{2\pi}{\displaystyle \int_{-\infty}^{\infty}\dfrac{e^{i\omega z}\left[\Gamma\left(1-i\omega\right)\right]^{\left(n+1\right)/2}\left[\Gamma\left(1+i\omega\right)\right]^{\left(n-1\right)/2}}{i\omega}d\omega}
\]
\[
+\left(-\pi i\right)\textrm{Res}\left[\dfrac{e^{i\omega z}\left[\Gamma\left(1-i\omega\right)\right]^{\left(n+1\right)/2}\left[\Gamma\left(1+\omega\right)\right]^{\left(n-1\right)/2}}{2\pi i\omega}\right]_{\omega=0}
\]
\[
=\dfrac{1}{2\pi}{\displaystyle \int_{-\infty}^{\infty}\dfrac{e^{i\omega z}\left[\Gamma\left(1-i\omega\right)\right]^{\left(n+1\right)/2}\left[\Gamma\left(1+i\omega\right)\right]^{\left(n-1\right)/2}}{i\omega}d\omega}.
\]
Finally, making the substitution $\varpi=i\omega$ gives
\[
F_{Z_{n}}\left(z\right)=\dfrac{1}{2\pi i}{\displaystyle \int_{-\infty i}^{\infty i}\dfrac{e^{\varpi z}\left[\Gamma\left(1-\varpi\right)\right]^{\left(n+1\right)/2}\left[\Gamma\left(1+\varpi\right)\right]^{\left(n-1\right)/2}}{\varpi}d\varpi}
\]
\[
=\dfrac{1}{2\pi i}{\displaystyle \int_{-\infty i}^{\infty i}e^{\varpi z}\left[\Gamma\left(1-\varpi\right)\right]^{\left(n+1\right)/2}\Gamma\left(1-1+\varpi\right)\left[\Gamma\left(1-0+\varpi\right)\right]^{\left(n-3\right)/2}d\varpi},
\]
which is equivalent to the first two lines of (20).

For part (ii), we present the derivation for $n\in\left\{ 1,3,\ldots\right\} $,
again noting that the simpler case of $n\in\left\{ 2,4,\ldots\right\} $
is analogous. To this end, insert the first line of (18) into (A1),
yielding
\[
F_{Z_{n}}\left(z\right)=\dfrac{1}{2}+\dfrac{1}{2\pi}{\displaystyle \int_{0}^{\infty}\left(\dfrac{\pi\omega}{\sinh\left(\pi\omega\right)}\right)^{\left(n-1\right)/2}\left\{ \dfrac{e^{i\omega z}{\displaystyle {\textstyle \int_{0}^{\infty}}e^{-t}\left[\cos\left(\ln\left(t\right)\omega\right)-i\sin\left(\ln\left(t\right)\omega\right)\right]dt}}{i\omega}\right.}
\]
\[
\left.-\dfrac{e^{-i\omega z}{\displaystyle {\textstyle \int_{0}^{\infty}}e^{-t}\left[\cos\left(\ln\left(t\right)\omega\right)+i\sin\left(\ln\left(t\right)\omega\right)\right]dt}}{i\omega}\right\} d\omega
\]
\[
=\dfrac{1}{2}+\dfrac{1}{2\pi}{\displaystyle \int_{0}^{\infty}\left(\dfrac{\pi\omega}{\sinh\left(\pi\omega\right)}\right)^{\left(n-1\right)/2}\left\{ \dfrac{\left[\cos\left(z\omega\right)+i\sin\left(z\omega\right)\right]{\displaystyle {\textstyle \int_{0}^{\infty}}e^{-t}\left[\cos\left(\ln\left(t\right)\omega\right)-i\sin\left(\ln\left(t\right)\omega\right)\right]dt}}{i\omega}\right.}
\]
\[
\left.-\dfrac{\left[\cos\left(z\omega\right)-i\sin\left(z\omega\right)\right]{\displaystyle {\textstyle \int_{0}^{\infty}}e^{-t}\left[\cos\left(\ln\left(t\right)\omega\right)+i\sin\left(\ln\left(t\right)\omega\right)\right]dt}}{i\omega}\right\} d\omega
\]
\[
=\dfrac{1}{2}+\dfrac{1}{2\pi}{\displaystyle \int_{0}^{\infty}\left(\dfrac{\pi\omega}{\sinh\left(\pi\omega\right)}\right)^{\left(n-1\right)/2}\left\{ \dfrac{{\displaystyle {\textstyle \int_{0}^{\infty}}e^{-t}\left[\cos\left(\left[z-\ln\left(t\right)\right]\omega\right)+i\sin\left(\left[z-\ln\left(t\right)\right]\omega\right)\right]dt}}{i\omega}\right.}
\]
\[
\left.-\dfrac{{\displaystyle {\textstyle \int_{0}^{\infty}}e^{-t}\left[\cos\left(\left[z-\ln\left(t\right)\right]\omega\right)-i\sin\left(\left[z-\ln\left(t\right)\right]\omega\right)\right]dt}}{i\omega}\right\} d\omega
\]
\[
=\dfrac{1}{2}+\dfrac{1}{\pi}{\displaystyle \int_{0}^{\infty}\left(\dfrac{\pi\omega}{\sinh\left(\pi\omega\right)}\right)^{\left(n-1\right)/2}\dfrac{{\displaystyle {\textstyle \int_{0}^{\infty}}e^{-t}\sin\left(\left[z-\ln\left(t\right)\right]\omega\right)dt}}{\omega}}d\omega.
\]
After rearrangement, this is equivalent to the first line of (21).
$\blacksquare$
\end{document}